\documentclass{article}
\usepackage{amsmath}
\usepackage{amssymb}
\usepackage{amsfonts}

\begin{document}

\title{Non-commutative probability and non-commutative processes: beyond the
Heisenberg algebra}
\author{R. Vilela Mendes\thanks{%
rvilela.mendes@gmail.com; rvmendes@fc.ul.pt} \\
CMAFCIO and IPFN, Universidade de Lisboa}
\date{ }
\maketitle

\begin{abstract}
A probability space is a pair ($\mathcal{A},\phi $) where $\mathcal{A}$ is
an algebra and $\phi $ a state on the algebra. In classical probability $%
\mathcal{A}$ is the algebra of linear combinations of indicator functions on
the sample space and in quantum probability $\mathcal{A}$ is the Heisenberg
or Clifford algebra. However, other algebras are of interest in
non-commutative probability. After a short review of the framework of
classical and quantum probability, other non-commutative probability spaces
are discussed, in particular those associated to non-commutative space-time.
\end{abstract}

\section{Introduction}

Probability plays a prominent role in physics and in all other natural
sciences. This is significant not only for the working scientist but also
for foundational questions. Most important issues are the interpretation of
probability results as well as the existence of one or more probability
frameworks. Quantum mechanics has called our attention to the need to
generalize the classical probability framework, the main structuring role
there being played by the Heisenberg algebra representations. This clearly
suggests that other probability structures might be associated to other
algebras relevant to phenomena in the natural sciences. In this paper, after
a short introduction to the formulation of noncommutative probability and
the particular case of the (Heisenberg algebra) quantum probability, a
construction is made of probability structures associated to $iso(1,1)$ and $%
iso(3,1)$ algebras, which arise in some formulations of noncommutative
spacetime.

\subsection{Non-commutative probability}

In classical probability theory, a \textit{probability space} is a triple $%
\left( \Omega ,F,P\right) $ where $\Omega $ (the sample space) is the set of
all possible outcomes, $F$ (the set of events) a $\sigma -$algebra of
subsets of $\Omega $ and $P$ a countably additive function from $F$ to $%
\left[ 0,1\right] $ assigning probabilities to events.

Let $H\mathcal{=}L^{2}\left( \Omega ,P\right) $ be the Hilbert space of
square-integrable functions on $\Omega $. In $H$ consider the set of \textit{%
indicator functions} $X_{B}=I_{\left\{ x\in B\right\} }$, with $x\in \Omega $
and $B\in F$. Given an unit vector \textbf{1} in $H$, the \textit{law} of $%
X_{B}$ is the probability measure%
\begin{equation}
B\rightarrow \left( 1,X_{B}1\right) =\int X_{B}\left( x\right) dP\left(
x\right) =P\left( B\right) .  \label{1.1}
\end{equation}%
The $X_{B}$'s form a set of \textit{countably additive projections} on
functions that only depend on the events $B$. $\left\{ X_{B}\right\} $ with 
\textit{involution} (by complex conjugation) is a $\ast -$\textit{algebra} $%
\mathcal{A}$, the algebra of random variables. The algebra may be identified
(by homomorphism) with the $\ast -$algebra $\mathcal{L}\left( H\right) $ of
bounded operators in $H$.

A \textit{spectral measure} $X$ over a measurable space $E$ is a countably
additive mapping from measurable sets to projections in some Hilbert space.
In this way classical probability is reframed as a spectral measure over $%
\Omega $. When a non-negative real number is assigned to each positive
element in the algebra $\mathcal{A}$ (such functional is called a \textit{%
state} in $\mathcal{A}$) one obtains von Neuman's view of probability theory.

In classical probability theory $\mathcal{A}$ is a commutative algebra.
Replacing it by a non-commutative $\ast -$algebra and keeping the
prescription of obtaining the law of $X\in \mathcal{A}$ by computation in a
state as in (\ref{1.1}) one obtains \textit{non-commutative probability}.

\subsection{Non-commutative processes}

Going from probability spaces to stochastic processes a few more steps are
required. Let us revisit Kolmogorov's extension theorem. It states that
given some interval $T$ in the real line, a finite sequence of points $%
t_{1},t_{2},\cdots ,t_{n}$ in this interval, a probability measure $\nu
_{t_{1},\cdots ,t_{n}}$ in $\mathbb{R}^{n}$ and measurable sets $%
B_{1},\cdots ,B_{n}$ satisfying

1 - $\nu _{\pi \left( t_{1}\right) ,\cdots ,\pi \left( t_{n}\right) }\left(
B_{\pi \left( t_{1}\right) }\times \cdots \times B_{\pi \left( t_{n}\right)
}\right) =\nu _{t_{1},\cdots ,t_{n}}\left( B_{t_{1}}\times \cdots \times
B_{t_{n}}\right) $ for any permutation $\pi $

2 - $\nu _{t_{1},\cdots ,t_{n}}\left( B_{t_{1}}\times \cdots \times
B_{t_{n}}\right) =\nu _{t_{1},\cdots ,t_{n},t_{n+1},\cdots ,t_{n+k}}\left(
B_{t_{1}}\times \cdots \times B_{t_{n}}\times \underset{k}{\underbrace{%
\mathbb{R}\cdots \times \mathbb{R}}}\right) $, then a probability space $%
\left( \Omega ,\mathcal{F},P\right) $ exists, as well as a stochastic
process $X_{t}:T\times \Omega \rightarrow \mathbb{R}$ such that%
\begin{equation*}
\nu _{t_{1},\cdots ,t_{n}}\left( B_{t_{1}}\times \cdots \times
B_{t_{n}}\right) =P\left( X_{t_{1}}\in B_{1},\cdots ,X_{t_{n}}\in
B_{n}\right) .
\end{equation*}

The index set $T$ may be thought of a part of the spectrum of an operator.
In the commutative case it is a subset of a multiplicative operator $%
\widehat{T}$. When going to non-commutative probability if the $\widehat{T}$
operator commutes with the elements of the probability algebra $\mathcal{A}$
then the construction of the non-commutative process proceeds in a way
similar to the Kolmogorov extension theorem. For each $t_{i}$ a
non-commutative probability space $\left( \mathcal{A}_{i},P_{i}\right) $, $%
P_{i}$ being a state, is constructed and then compatibility conditions as in
the Kolmogorov's theorem are imposed. Processes of this type, of which 
\textit{quantum probability} is an example, will be called \textit{%
non-commutative processes of type I}.

If however the index set operator $\widehat{T}$ has nontrivial commutation
relations with the elements of the probability algebra $\mathcal{A}$, the
construction of the process will be different. Processes where $\widehat{T}$
does not commute with $\mathcal{A}$ will be called \textit{non-commutative
processes of type II}.

In the next sections, after a short review of the quantum probability
setting, non-commutative processes for more general algebras will be
discussed as well as the construction of non-commutative processes of type
II. Some results in this context were reported in \cite{Vilela3}. Here a
more systematic presentation is given as well as some new results.

\section{Quantum probability: Non-commutative processes and the Heisenberg
algebra}

The main motivation to extend probability theory to the non-commutative
setting came from the application
of~probabilistic~concepts~such~as~independence and~noise to quantum
mechanics. For this reason the class of developments in non-commutative
probability inspired by the structure of Quantum Mechanics carries the names
of \textit{Quantum Probability} or \textit{Quantum Stochastic Processes} or,
more generally, \textit{Quantum Stochastic Analysis} \cite{Accardi1} \cite%
{Partha1} \cite{Partha2} \cite{Meyer1}.

The dynamics of a particle in classical mechanics is described in phase
space by functions of its coordinate $q$ and momentum $p$. Functions $%
f\left( q,p\right) $ form a commutative algebra and a classical
probabilistic description of the particle dynamics is an assignment of a
probability $P_{f}$ to each function $f\left( q,p\right) $ by%
\begin{equation}
P_{f}=\int f\left( p,q\right) d\mu \left( p,q\right) ,  \label{2.1}
\end{equation}%
the measure $\mu \left( p,q\right) $, with $\int d\mu \left( p,q\right) =1$,
being the state.

In Quantum Mechanics such functions do not commute. In particular, for the
phase-space coordinate functions one has the canonical commutation relations
(CCR),%
\begin{equation}
\left[ q,p\right] =i\hslash  \label{2.2}
\end{equation}%
with $\hslash =1.054\times 10^{-34}J\times s=1.054\times 10^{-27}g\times
cm^{2}\times s^{-1}$. For phenomena at the scale of $cm$ ($q$) and $g\times
cm\times s^{-1}$ ($p$), $\hslash $\ is a small quantity. Nevertheless, being
non zero, it entirely changes the structure. For $n$ particle species it
would be%
\begin{eqnarray}
\left[ q_{i},p_{j}\right] &=&i\hslash \delta _{ij}\mathbb{I}  \notag \\
\left[ q_{i},q_{j}\right] &=&\left[ p_{i},p_{j}\right] =\left[ p_{i},\mathbb{%
I}\right] =\left[ x_{i},\mathbb{I}\right] =0.  \label{2.3}
\end{eqnarray}%
This is the Lie algebra of the Heisenberg group $\mathcal{H}\left( n\right) $%
, the maximal nilpotent subgroup in $U\left( n+1,1\right) $. \textit{Quantum
Probability} is the particular case of non-commutative probability
associated to this Heisenberg algebra. Let $n=1$. In the unitary
representations of $\mathcal{H}\left( 1\right) $ in $L^{2}\left( \mathbb{R}%
\right) $ the Lie algebra operators are (\cite{Vilenkin} Ch.12)%
\begin{equation}
\begin{split}
qf\left( x\right) =& xf\left( x\right) \\
pf\left( x\right) =& -i\lambda \frac{d}{dx}f\left( x\right) \\
\mathbb{I}f\left( x\right) =& \lambda f\left( x\right)
\end{split}
\label{2.4}
\end{equation}%
the representations being irreducible for $\lambda \neq 0$. Defining
creation and annihilation operators%
\begin{eqnarray}
a &=&\frac{1}{\sqrt{2}}\left( q+ip\right)  \notag \\
a^{\dag } &=&\frac{1}{\sqrt{2}}\left( q-ip\right)  \label{2.5}
\end{eqnarray}%
one obtains a representation of the operators in Fock space, a convenient
dense set in Fock space being the set of exponential vectors%
\begin{equation}
\psi \left( f\right) =1\oplus f\oplus \cdots \oplus \frac{f^{(n)}}{\sqrt{n!}}%
\oplus \cdots  \label{2.6}
\end{equation}%
where $f^{(n)}$ is the $n-$fold tensor product of $f$.

To define processes one considers an index set $\left[ 0,T\right) $ and a
family of operators $\left\{ A_{t},A_{t}^{\dag },\Lambda _{t}\right\} $
indexed by the characteristic functions $\left[ 0,t\right) $. These
operators have the following action on exponential vectors.%
\begin{eqnarray}
A_{t}\psi \left( f\right) &=&\int_{0}^{t}dsf\left( s\right) \psi \left(
f\right)  \notag \\
A_{t}^{\dag }\psi \left( f\right) &=&\left. \frac{\partial }{\partial
\varepsilon }\right\vert _{\varepsilon =0}\psi \left( f+\varepsilon \chi _{%
\left[ 0,t\right] }\right)  \notag \\
\Lambda _{t}\psi \left( f\right) &=&\left. \frac{\partial }{\partial
\varepsilon }\right\vert _{\varepsilon =0}\psi \left( e^{\varepsilon \chi _{%
\left[ 0,t\right] }}f\right) .  \label{2.7}
\end{eqnarray}%
Quantum stochastic differentials are defined by%
\begin{eqnarray*}
dA_{t} &=&A_{t+dt}-A_{t} \\
dA_{t}^{\dag } &=&A_{t+dt}^{\dag }-A_{t}^{\dag } \\
d\Lambda _{t} &=&\Lambda _{t+dt}-\Lambda _{t}
\end{eqnarray*}%
and quantum stochastic calculus is, in practice, an application of the rules%
\begin{equation*}
dA_{t}dA_{t}^{\dag }=dt;\;dA_{t}d\Lambda _{t}=dA_{t};\;d\Lambda
_{t}dA_{t}^{\dag }=dA_{t}^{\dag };\;d\Lambda _{t}d\Lambda _{t}=d\Lambda _{t}
\end{equation*}%
all other products vanishing.

Many deep results have been obtained in the quantum stochastic processes
field\footnote{%
See for example the proceedings of the serie of conferences "Quantum
Probability and Infinite Dimensional Analysis" and references therein.} of
practical importance for physical systems perturbed by quantum noise. Also,
nonlinear extensions have been obtained \cite{Accardi2} \cite{Accardi3} \cite%
{Accardi4}.

In quantum mechanics bosons satisfy CCR whereas fermions have canonical
anticommutation relations (CAR). Based on the fermion Clifford algebra a
non-commutative stochastic calculus has also been developed \cite{Streater1} 
\cite{Streater2} \cite{Carlen}. Nevertheless boson and fermion stochastic
calculus may be unified in a single theory \cite{Hudson}.

Other extensions, closely related to the Heisenberg algebra, have been made
by Boukas \cite{Boukas} who deals with a discrete analogue of the Heisenberg
algebra and Privault \cite{Privault} for L\'{e}vy processes on finite
difference algebras. Generalized Gaussian processes have also been
constructed associated to the q-deformed Fock space \cite{Bozejko1} or to
extensions of this concept \cite{Bozejko2} \cite{Ejsmont} \cite{Bozejko3}.

\section{Non-commutative processes beyond the Heisenberg and Clifford
algebras}

\subsection{Non-commutative processes for the iso(1,1) algebra}

This probability space was first discussed in \cite{Vilela3} and the main
results are recalled here. The algebra is%
\begin{eqnarray}
\left[ P,Q\right] &=&-i\mathfrak{I}  \notag \\
\left[ Q,\mathfrak{I}\right] &=&-iP  \notag \\
\left[ P,\mathfrak{I}\right] &=&0.  \label{3.1}
\end{eqnarray}%
This is the algebra of one-dimensional non-commutative phase-space with the
correspondence to the physical variables $p,x$ established by $P=\mathfrak{%
\ell }p$ and $Q=\frac{x}{\mathfrak{\ell }}$, $\mathfrak{\ell }$ being a
fundamental length.

As in quantum probability, the first step is to obtain the representations
of the algebra. The algebra (\ref{3.1}) is the Lie algebra of the group of
motions of pseudo-Euclidean space $ISO(1,1)$. The group contains hyperbolic
rotations in the plane and two translations, being isomorphic to the group
of $3\times 3$ matrices%
\begin{equation}
\left( 
\begin{array}{ccc}
\cosh \mu & \sinh \mu & a_{1} \\ 
\sinh \mu & \cosh \mu & a_{2} \\ 
0 & 0 & 1%
\end{array}%
\right)  \label{3.2}
\end{equation}%
its elements being labelled by $\left( \mu ,a_{1},a_{2}\right) $. The
representations of the group in the space of functions in the hyperbola $%
f\left( \cosh \nu ,\sinh \nu \right) $ are%
\begin{equation}
T_{R}\left( \mu ,a_{1},a_{2}\right) f\left( \nu \right) =e^{R\left(
-a_{1}\cosh \mu +a_{2}\sinh \mu \right) }f\left( \nu -\mu \right) .
\label{3.3}
\end{equation}%
The $T_{R}$ representation is unitary if $R=ir$ is purely imaginary and
irreducible if $R\neq 0$. For the Lie algebra one obtains:%
\begin{eqnarray}
Q &=&i\frac{d}{d\mu }  \notag \\
P &=&\sinh \mu  \notag \\
\mathfrak{I} &\mathfrak{=}&\cosh \mathfrak{\mu }  \label{3.4}
\end{eqnarray}%
acting on the space $V_{1}$ of functions on the hyperbola. The $r=1$
representation has been chosen. A quite similar construction applies to all $%
r\neq 0$ cases. It is convenient to define the operators%
\begin{eqnarray}
A_{+} &=&\frac{1}{\sqrt{2}}\left( Q-iP\right)  \notag \\
A_{-} &=&\frac{1}{\sqrt{2}}\left( Q+iP\right) .  \label{3.5}
\end{eqnarray}%
Using the representation (\ref{3.4}) the state in $V_{1}$ that is
annihilated by $A_{-}$ is%
\begin{equation}
\psi _{0}=\frac{1}{\sqrt{N}}e^{-\cosh \mu }  \label{3.6}
\end{equation}%
$N$ being a normalization factor $N=2K_{0}\left( 2\right) $, $K_{0}$ being a
modified Bessel function of the 2nd kind.

With the algebra of the operators $\left\{ A_{+},A_{-},\mathfrak{I}\right\} $
and the state defined by $\psi _{0}-$expectation values the probability
space is defined. The rest of the construction is standard.

With the scalar product%
\begin{equation}
\left( \psi ,\phi \right) _{\mu }=\int_{R}d\mu \psi ^{\ast }\left( \mu
\right) \phi \left( \mu \right)  \label{3.7}
\end{equation}%
the space $V_{1}$ of square integrable functions in the hyperbola becomes,
by completion, a Hilbert space and let $h=L^{2}\left( R_{+}\right) $ be the
Hilbert space of square integrable functions on the half-line $R_{+}=\left[
0,\infty \right) $. One now constructs an infinite set of operators labelled
by functions on $h$, with the algebra%
\begin{eqnarray}
\left[ A_{-}\left( f\right) ,A_{+}\left( g\right) \right] &=&\mathfrak{I}%
\left( fg\right)  \notag \\
\left[ A_{-}\left( f\right) ,\mathfrak{I}\left( g\right) \right] &=&\left[
A_{+}\left( f\right) ,\mathfrak{I}\left( g\right) \right] =\frac{1}{2}\left(
A_{-}\left( fg\right) -A_{+}\left( fg\right) \right) .  \notag
\end{eqnarray}%
These operators act on $H_{1}=h\otimes V_{1}$ by%
\begin{eqnarray}
A_{-}\left( f\right) \psi \left( g\right) &=&\frac{i}{\sqrt{2}}\left( \left( 
\frac{d}{d\mu }+\sinh \mu \right) \psi \right) \left( fg\right)  \notag \\
A_{+}\left( f\right) \psi \left( g\right) &=&\frac{i}{\sqrt{2}}\left( \left( 
\frac{d}{d\mu }-\sinh \mu \right) \psi \right) \left( fg\right)  \notag \\
\mathfrak{I}\left( f\right) \psi \left( g\right) &=&\left( \cosh \mu \times
\psi \right) \left( fg\right)  \label{3.9}
\end{eqnarray}%
$f,g\in h$, $\psi \in V_{1}$, $\psi \left( g\right) \in $ $H_{1}$.

Adapted processes are associated to the splitting%
\begin{equation}
h=L^{2}\left( 0,t\right) \oplus L^{2}\left( t,\infty \right) =h^{t}\oplus
h^{(t}  \label{3.10}
\end{equation}%
with the corresponding%
\begin{equation}
H_{1}=H_{1}^{t}\oplus H_{1}^{(t}=h^{t}\otimes V_{1}\oplus h^{(t}\otimes
V_{1}.  \label{3.11}
\end{equation}%
An adapted process is a family $O=\left( O\left( t\right) ,t\geq 0\right) $
of operators such that $O\left( t\right) =O^{t}\otimes 1$. The basic adapted
processes are%
\begin{eqnarray}
\mathfrak{I}\left( t\right) &=&\mathfrak{I}\left( \chi _{\left[ 0,t\right]
}\right)  \notag \\
A_{+}\left( t\right) &=&A_{+}\left( \chi _{\left[ 0,t\right] }\right)  \notag
\\
A_{-}\left( t\right) &=&A_{-}\left( \chi _{\left[ 0,t\right] }\right)
\label{3.12}
\end{eqnarray}%
For each adapted process $O\left( t\right) $ define $dO\left( t\right)
=O\left( t+dt\right) -O\left( t\right) $. Given adapted processes $%
E_{0},E_{+},E_{-}$ a stochastic integral%
\begin{equation}
N\left( t\right) =\int_{0}^{t}E_{0}d\mathfrak{I}+E_{+}dA_{+}+E_{-}dA_{-}
\label{3.13}
\end{equation}%
is defined as the limit of the family of operators $\left( N\left( s\right)
;s\geq 0\right) $ such that $N\left( 0\right) =0$ and for $t_{n}<t\leq
t_{n}+dt$%
\begin{equation}
N\left( t\right) =N\left( t_{n}\right) +E_{0}^{(n)}\left( \mathfrak{I}\left(
t\right) -\mathfrak{I}\left( t_{n}\right) \right) +E_{+}^{(n)}\left(
A_{+}\left( t\right) -A_{+}\left( t_{n}\right) \right) +E_{-}^{(n)}\left(
A_{-}\left( t\right) -A_{-}\left( t_{n}\right) \right)  \label{3.15}
\end{equation}%
with $E_{0},E_{+},E_{-}$ assumed to be simple processes%
\begin{eqnarray}
E_{0} &=&\sum_{n=0}^{\infty }E_{0}^{(n)}\chi _{\left[ t_{n},t_{n+1}\right) }
\notag \\
E_{+} &=&\sum_{n=0}^{\infty }E_{+}^{(n)}\chi _{\left[ t_{n},t_{n+1}\right) }
\notag \\
E_{-} &=&\sum_{n=0}^{\infty }E_{-}^{(n)}\chi _{\left[ t_{n},t_{n+1}\right) }
\label{3.16}
\end{eqnarray}%
Then, from (\ref{3.15}) and (\ref{3.16}) it follows that in $H_{1}=h\otimes
V_{1}$ one has%
\begin{eqnarray*}
\left\langle \phi \left( f\right) ,N\left( t\right) \psi \left( g\right)
\right\rangle &=&\int_{0}^{t}\left\langle \phi \left( f\right) ,E_{0}\left(
s\right) \left( \cosh \mu \times \psi \right) \left( g\chi _{ds}\right)
\right. \\
&&+E_{-}\left( s\right) \frac{i}{\sqrt{2}}\left( \left( \frac{d}{d\mu }%
+\sinh \mu \right) \psi \right) \left( g\chi _{ds}\right) \\
&&\left. +E_{+}\left( s\right) \frac{i}{\sqrt{2}}\left( \left( \frac{d}{d\mu 
}-\sinh \mu \right) \psi \right) \left( g\chi _{ds}\right) \right\rangle .
\end{eqnarray*}%
Multiplication rules for the stochastic differentials 
\begin{eqnarray}
dA_{-}\left( t\right) &=&A_{-}\left( t+dt\right) -A_{-}\left( t\right) 
\notag \\
dA_{+}\left( t\right) &=&A_{+}\left( t+dt\right) -A_{+}\left( t\right) 
\notag \\
d\mathfrak{I}\left( t\right) &=&\mathfrak{I}\left( t+dt\right) -\mathfrak{I}%
\left( t\right)  \label{3.18}
\end{eqnarray}%
are obtained taking into account the commutation relations and expectation
values in the $\psi _{0}$ state. Non-vanishing products are%
\begin{eqnarray}
dA_{-}\left( t\right) dA_{+}\left( t\right) &=&\mathfrak{I}\left( dt\right) 
\notag \\
d\mathfrak{I}\left( t\right) dA_{+}\left( t\right) &=&-dA_{-}\left( t\right)
d\mathfrak{I}\left( t\right) =-\frac{1}{2}\left( A_{-}\left( dt\right)
-A_{+}\left( dt\right) \right) .  \label{3.19}
\end{eqnarray}

Summarizing, this proves the following:

\textbf{Proposition 1: }\textit{Associated to each irreducible unitary
representation of the }$iso(1,1)$\textit{\ algebra there is a stochastic
process characterized by the operator set }$\left\{ A_{+}\left( f\right)
,A_{-}\left( f\right) ,\mathfrak{I}\left( g\right) \right\} $\textit{\
acting on }$H_{1}=h\otimes V_{1}$\textit{\ with adapted processes and
stochastic integrals satisfying properties (\ref{3.12}-\ref{3.19}).}

Of interest is also the characterization of the coordinate process $Q\left(
t\right) =\frac{1}{\sqrt{2}}\left( A_{-}\left( t\right) +A_{+}\left(
t\right) \right) $. The characteristic functional is%
\begin{equation}
C_{Q\left( t\right) }\left( f\right) =\frac{K_{0}\left( 2\cos \left( \frac{1%
}{2}\int_{0}^{t}f\left( s\right) ds\right) \right) }{K_{0}\left( 2\right) }
\label{3.21}
\end{equation}

\subsection{Non-commutative processes for the iso(3,1) algebra}

The algebra here is%
\begin{eqnarray}
\left[ Q^{i},Q^{j}\right] &=&-iM^{ij}  \notag \\
\left[ P^{i},Q^{j}\right] &=&-i\mathfrak{I\delta }^{ij}  \notag \\
\left[ Q^{i},\mathfrak{I}\right] &=&iP^{i}  \notag \\
\left[ M^{ij},Q^{k}\right] &=&i\left( Q^{i}\delta ^{jk}-Q^{j}\delta
^{ik}\right)  \notag \\
\left[ M^{ij},P^{k}\right] &=&i\left( P^{i}\delta ^{jk}-P^{j}\delta
^{ik}\right)  \notag \\
\left[ P^{i},\mathfrak{I}\right] &=&\left[ M^{ij},\mathfrak{I}\right] =\left[
P^{i},P^{j}\right] =0.  \label{4.1}
\end{eqnarray}%
This is, in three dimensions, the phase space algebra of non-commutative
space-time. It is the algebra of the Poincar\'{e} group. For those familiar
with the use of this group in relativistic quantum mechanics notice that
here the three $Q^{i}$'s play the role of the boost generators and $%
\mathfrak{I}$ the role of the generator of time translations.

The group $ISO(3,1)$ is isomorphic to the group of matrices $\left( 
\begin{array}{cc}
M & a \\ 
0 & 1%
\end{array}%
\right) $ where $M\in SO(3,1)$ and $a\in R^{4}$.

A representation of the $iso\left( 3,1\right) $ algebra that generalizes the
representation (\ref{3.4}) is obtained in the space $V_{3}$ of functions in
the upper sheet of the hyperboloid $\mathcal{H}_{+}^{3}$ with coordinates%
\begin{eqnarray}
\xi _{1} &=&\sinh \mu \sin \theta _{1}\cos \theta _{2}  \notag \\
\xi _{2} &=&\sinh \mu \sin \theta _{1}\sin \theta _{2}  \notag \\
\xi _{3} &=&\sinh \mu \cos \theta _{1}  \notag \\
\xi _{4} &=&\cosh \mu .  \label{4.2}
\end{eqnarray}%
The representation is%
\begin{eqnarray}
\mathfrak{I} &\mathfrak{=}&\cosh \mathfrak{\mu }  \notag \\
P^{1} &=&\sinh \mu \sin \theta _{1}\cos \theta _{2}  \notag \\
P^{2} &=&\sinh \mu \sin \theta _{1}\sin \theta _{2}  \notag \\
P^{3} &=&\sinh \mu \cos \theta _{1}  \notag \\
M^{12} &=&i\frac{\partial }{\partial \theta _{2}}  \notag \\
M^{31} &=&i\left( \cos \theta _{2}\frac{\partial }{\partial \theta _{1}}%
-\cot \theta _{1}\sin \theta _{2}\frac{\partial }{\partial \theta _{2}}%
\right)  \notag \\
M^{23} &=&-i\left( \sin \theta _{2}\frac{\partial }{\partial \theta _{1}}%
+\cot \theta _{1}\cos \theta _{2}\frac{\partial }{\partial \theta _{2}}%
\right)  \notag \\
Q^{1} &=&i\left( \sin \theta _{1}\cos \theta _{2}\frac{\partial }{\partial
\mu }+\coth \mu \left( \cos \theta _{1}\cos \theta _{2}\frac{\partial }{%
\partial \theta _{1}}-\frac{\sin \theta _{2}}{\sin \theta _{1}}\frac{%
\partial }{\partial \theta _{2}}\right) \right)  \notag \\
Q^{2} &=&i\left( \sin \theta _{1}\sin \theta _{2}\frac{\partial }{\partial
\mu }+\coth \mu \left( \cos \theta _{1}\sin \theta _{2}\frac{\partial }{%
\partial \theta _{1}}+\frac{\cos \theta _{2}}{\sin \theta _{1}}\frac{%
\partial }{\partial \theta _{2}}\right) \right)  \notag \\
Q^{3} &=&i\left( \cos \theta _{1}\frac{\partial }{\partial \mu }-\coth \mu
\sin \theta _{1}\frac{\partial }{\partial \theta _{1}}\right)  \label{4.3}
\end{eqnarray}%
Define the operators%
\begin{eqnarray}
A_{-}^{k} &=&\frac{1}{\sqrt{2}}\left( Q^{k}+iP^{k}\right)  \notag \\
A_{+}^{k} &=&\frac{1}{\sqrt{2}}\left( Q^{k}-iP^{k}\right)  \label{4.4}
\end{eqnarray}%
to which in $H_{3}=h\otimes V_{3}$ correspond the operators $A_{-}^{k}\left(
t\right) \doteq A_{-}^{k}\left( \chi _{\left[ 0,t\right] }\right)
,A_{+}^{k}\left( t\right) \doteq A_{-}^{k}\left( \chi _{\left[ 0,t\right]
}\right) $ as well as $\mathfrak{I}\left( t\right) \mathfrak{\doteq I}\left(
\chi _{\left[ 0,t\right] }\right) ,M^{kj}\left( t\right) \doteq M^{kj}\left(
\chi _{\left[ 0,t\right] }\right) $ with algebra 
\begin{eqnarray}
\left[ A_{-}^{k}\left( t\right) ,A_{+}^{j}\left( s\right) \right] &=&\delta
^{kj}\mathfrak{I}\left( t\wedge s\right) -\frac{i}{2}M^{kj}\left( t\wedge
s\right)  \notag \\
\left[ A_{-}^{k}\left( t\right) ,\mathfrak{I}\left( s\right) \right] &=&%
\left[ A_{+}^{k}\left( t\right) ,\mathfrak{I}\left( s\right) \right] =\frac{i%
}{\sqrt{2}}P^{k}\left( t\wedge s\right)  \notag \\
\left[ A_{\pm }^{k}\left( t\right) ,M^{jn}\left( s\right) \right] &=&i\left(
\delta ^{kj}A_{\pm }^{n}\left( t\wedge s\right) -\delta ^{kn}A_{\pm
}^{j}\left( t\wedge s\right) \right)  \notag \\
\left[ M^{kj}\left( t\right) ,M^{nl}\left( s\right) \right] &=&-i\left(
\delta ^{jn}M^{kl}+\delta ^{kl}M^{jn}-\delta ^{kn}M^{jl}-\delta
^{jl}M^{kn}\right) \left( t\wedge s\right) ,  \notag \\
&&  \label{4.5}
\end{eqnarray}%
which follows from the commutation relations (\ref{4.1}). The normalized
spherical symmetrical state $\psi _{0}=\frac{1}{\sqrt{8\pi K_{0}\left(
2\right) }}\exp \left( -\cosh \mu \right) $ is annihilated by all $%
A_{-}^{k}\left( t\right) $%
\begin{equation}
A_{-}^{k}\left( t\right) \psi _{0}=0;\;k=1,2,3  \label{4.6}
\end{equation}%
In conclusion:

\textbf{Proposition 2: }\textit{The operators }$\left\{ A_{-}^{k}\left(
t\right) ,A_{+}^{j}\left( t\right) ,\mathfrak{I}\left( t\right)
,M^{kj}\left( t\right) \right\} $\textit{, the representation (\ref{4.3})
and expectations on the }$\psi _{0}$\textit{\ state define a probability
space. Multiplication rules for the stochastic differentials }$%
dA_{-}^{k}\left( t\right) =A_{-}^{k}\left( t+dt\right) -A_{-}^{k}\left(
t\right) $\textit{, etc. are}%
\begin{eqnarray}
dA_{-}^{k}\left( t\right) dA_{+}^{j}\left( t\right) &=&\delta ^{kj}\mathfrak{%
I}\left( dt\right)  \notag \\
dA_{+}^{k}\left( t\right) dA_{-}^{j}\left( t\right) &=&0  \notag \\
dA_{-}^{k}\left( t\right) d\mathfrak{I}\left( t\right) &=&-d\mathfrak{I}%
\left( t\right) dA_{-}^{k}\left( t\right) =\frac{1}{2}\left( A_{-}^{k}\left(
dt\right) -A_{+}^{k}\left( dt\right) \right) .  \label{4.7}
\end{eqnarray}

The multiplication rules for the stochastic differentials are obtained, as
before, by taking into account the commutation relations and expectation
values in the $\psi _{0}$ state.

The characteristic functionals of the coordinate processes $Q^{k}\left(
t\right) =\frac{1}{\sqrt{2}}\left( A_{-}^{k}\left( t\right) +A_{+}^{k}\left(
t\right) \right) $ are%
\begin{equation}
C_{Q^{k}\left( t\right) }\left( f\right) =\frac{1}{4\pi K_{0}\left( 2\right) 
}\int d\Omega K_{0}\left( 2\cos \left( \frac{g_{k}\left( \theta _{1},\theta
_{2}\right) }{2}\int_{0}^{t}f\left( s\right) ds\right) \right)  \label{4.8}
\end{equation}%
with $g_{1}\left( \theta _{1},\theta _{2}\right) =\sin \theta _{1}\cos
\theta _{2};\;g_{2}\left( \theta _{1},\theta _{2}\right) =\sin \theta
_{1}\sin \theta _{2};\;g_{3}\left( \theta _{1},\theta _{2}\right) =\cos
\theta _{1}$ and $d\Omega =\sin \theta _{1}d\theta _{1}d\theta _{2}$. Notice
however that these processes are not independent being related by the
rotation process $M^{kj}\left( t\right) $%
\begin{equation}
\left[ Q^{k}\left( t\right) ,Q^{j}\left( s\right) \right] =-iM^{kj}\left(
t\wedge s\right) .  \label{4.9}
\end{equation}

\section{Non-commutative processes of type II}

In the previous section the index set of the stochastic processes (or the $%
h- $space) may be looked at as the continuous spectrum of an operator that
commutes with all other operators in the algebra. Here the case where such
operator does not commute with the algebra operators will be analyzed.

In this case the $H$ space, where the stochastic process operators act, can
no longer be a direct product as before, rather it is a representation space
for the whole set of operators, including the one that generates the index
set. An example is the following algebra%
\begin{eqnarray}
\left[ E,T\right] &=&i\mathfrak{I}  \notag \\
\left[ T,\mathfrak{I}\right] &=&iE  \notag \\
\left[ Q,P\right] &=&i\mathfrak{I}  \notag \\
\left[ Q,\mathfrak{I}\right] &=&iP  \notag \\
\left[ P,\mathfrak{I}\right] &=&\left[ P,T\right] =\left[ M^{10},\mathfrak{I}%
\right] =0  \notag \\
\left[ Q,T\right] &=&-iM^{10}  \notag \\
\left[ M^{10},Q\right] &=&iT  \notag \\
\left[ M^{10},P\right] &=&iE  \notag \\
\left[ M^{10},T\right] &=&iQ  \notag \\
\left[ M^{10},E\right] &=&iP.  \label{5.1}
\end{eqnarray}%
This $iso\left( 2,1\right) $ algebra is the algebra of two-dimensional
non-commutative space-time. Notice that as a real algebra it is different
from the algebra studied in \cite{Vilela3} because a different metric is
used in the deformed space. For consistency with the physical interpretation
of $T$ as the time operator, $T$ will be the operator that generates the
index set of the stochastic process. A representation of this algebra is
obtained on functions $f\left( \sigma ,\theta \right) $ in the $C^{2}$ cone,%
\begin{eqnarray}
T &=&-i\frac{\partial }{\partial \theta }  \notag \\
E &=&\sigma \sin \theta  \notag \\
\mathfrak{I} &=&\sigma \cos \theta  \notag \\
P &=&\sigma  \notag \\
Q &=&i\left( \sigma \cos \theta \frac{\partial }{\partial \sigma }-\sin
\theta \frac{\partial }{\partial \theta }\right)  \notag \\
M^{10} &=&i\left( \sigma \sin \theta \frac{\partial }{\partial \sigma }+\cos
\theta \frac{\partial }{\partial \theta }\right)  \label{5.2}
\end{eqnarray}%
$T$ has discrete spectrum, with eigenfunctions%
\begin{equation}
T\left( e^{in\theta }f\left( \sigma \right) \right) =ne^{in\theta }f\left(
\sigma \right)  \label{5.3}
\end{equation}%
$f\left( \sigma \right) $ being an arbitrary function of $\sigma $.

In the cases studied before the space $H$, where the stochastic process
operators act, is a direct product $h\otimes V$, that is, the space $V$ is
the same for all times. Here however once a normalized cyclic vector%
\footnote{$\exp \left( -\sigma ^{2}/2\right) /\left( \pi ^{1/4}\right) $ for
example} of $\sigma $ is chosen to generate a $v_{\sigma }$ space, the space 
$V_{n}$ at time $n$ is $V_{n}=e^{in\theta }v_{\sigma }$ and%
\begin{equation}
H=\oplus V_{n}  \label{5.4}
\end{equation}%
An adapted, discrete time, stochastic process is%
\begin{equation}
O\left( t\right) =\oplus _{n=0}^{t}O_{n}  \label{5.5}
\end{equation}%
$O_{n}$ being one of the operators in (\ref{5.2}) acting on $V_{n}$. The
fact that in this case the space $H$ is not a direct product is the main
difference from the previous cases. Otherwise the construction is similar
with (discrete) stochastic differentials defined by%
\begin{equation*}
dO\left( n\right) =O_{n+1}-O_{n}
\end{equation*}%
with $O_{n+1}$ operating in $V_{n+1}$ and $O_{n}$ on $V_{n}$.

For example, for functionals $\mathcal{F}\left( Q\left( t\right) \right) $
of the coordinate process $Q\left( t\right) $, expectation values are
obtained by%
\begin{equation}
\left\langle \mathcal{F}\left( Q\left( t\right) \right) \right\rangle
=\sum_{0}^{t}\frac{1}{2\pi }\int d\sigma d\theta \phi ^{\ast }\left( \sigma
\right) e^{-in\theta }\mathcal{F}\left( Q\left( \sigma ,\theta \right)
\right) e^{in\theta }\phi \left( \sigma \right) .  \label{5.6}
\end{equation}

\section{Remarks and conclusions}

1 - von Neuman's view of probability theory as a pair ($\mathcal{A},\phi $),
where $\mathcal{A}$ is an algebra and $\phi $ a state, is a powerful insight
with far reaching implications. The main purpose of this paper was to
emphasize that in addition to the many results already obtained in quantum
probability, where the Heisenberg or the Clifford algebras play the central
role, other algebras are of potential interest, in particular in the
non-commutative space-time context \cite{Vilela3} \cite{Vilela2} \cite%
{Vilela4}.

In addition, it has been pointed out that in some cases the index set of the
stochastic process might be generated by a non trivial operator in the
algebra, the consequence being that the space where the stochastic process
operators act is no longer a direct product.

Most results on the processes associated to the iso(1,1)-algebra were
previously reported in \cite{Vilela3} but not those related to the
iso(3,1)-algebra. Also, concerning the processes of type II a different
metric structure is used, which might be physically more relevant.

2 - In the past, a central role is played in Quantum Probability by the
algebra of the second quantization operators, the creation and annihilation
operators. In 1977 Palev \cite{Palev1} has shown that the basic requirements
of the second quantization formulation may also be obtained by operators
that lead to generalized quantum statistics \cite{Palev2}. The algebras of
these generalized quantum statistics operators would also provide new
examples of non-commutative processes.

\section*{Acknowledgements}
Supported by Funda\c{c}\~{a}o para a Ci\^{e}ncia e a Tecnologia, under the project: UID/MAT/04561/2019

\end{document}